\documentclass[12pt,a4paper]{article}
\usepackage{amsmath,amsfonts,amssymb,amscd}
\usepackage{xypic}

\usepackage[a4paper]{geometry}
\usepackage{indentfirst}

\usepackage{tikz}
\usetikzlibrary{matrix}

\def\H{\operatorname{H}}

\def\Ker{\operatorname{Ker}}
\def\Ad{\operatorname{Ad}}

\def\red{\operatorname{red}}
\def\d{\operatorname{d}}
\def\GL{\operatorname{GL}}

\def\OSp{\operatorname{OSp}}

\def\Im{\operatorname{Im}}

\newcounter{th}
\def\t{\refstepcounter{th}{\bf \noindent{Theorem} \arabic{th}. }}

\newcounter{prop}

\newcounter{lem}
\def\lem{\refstepcounter{lem}{\bf \noindent{Lemma} \arabic{lem}. }}

\newcounter{de}

\newcounter{ex}

\begin{document}

\begin{center}

{\LARGE{\bf Vector fields on $\mathfrak{osp}_{2m-1|2n}(\mathbb C)$-flag supermanifolds}}

\bigskip

{\bf Elizaveta Vishnyakova}\\[0.3cm]
\end{center}
\bigskip
\bigskip

\begin{abstract}
	We compute the Lie superalgebras of holomorphic vector fields on  isotropic flag supermanifolds  of maximal type corresponding to the Lie superalgebra $\mathfrak{osp}_{2m-1|2n}(\mathbb C)$. 
	
\end{abstract}

\bigskip

\section{Introduction}

 In \cite{Man} Yu.I.~Manin defined four series of complex homogeneous supermanifolds that correspond to four series of linear Lie superalgebras: $\mathfrak{gl}_{m|n}(\mathbb C)$, $\mathfrak{osp}_{m|2n}(\mathbb C)$, $\pi\mathfrak{sp}_{n}(\mathbb C)$ and $\mathfrak{q}_{n}(\mathbb C)$. The Lie superalgebras of holomorphic vector fields on these supermanifolds were studied in \cite{bunOni,buneg,oniq,onich,onigl,oniosp,onipisp,serov} in the case of super-Grassmannians and in \cite{ViGL in JA,ViPi-sym,Vivector,Vi osp and pi} in the case of other flag supermanifolds. We summarize the obtained results in the following table:
 \bigskip
  \begin{center}
 	\begin{tabular}{|c||c|c|}
 		\hline
 		& Super-Grassmannians &  Other
 		flag supermanifolds\\
 				\hline
 		$\mathfrak{gl}_{m|n}(\mathbb{C})$ & + & + (generic)  \\
 		\hline
 		$\mathfrak{osp}_{2m|2n}(\mathbb{C})$ & + (maximal type) &  + (maximal type)\\
 		\hline
 		$\mathfrak{osp}_{2m+1|2n}(\mathbb{C})$ & + (maximal type) &  - \\
 		\hline
 		$\pi\mathfrak{sp}_{n|n}(\mathbb{C})$ & + (maximal type) &  + (maximal type)\\
 		\hline
 		$\mathfrak{q}_{n|n}(\mathbb{C})$ & + & + \\
 		\hline
 	\end{tabular}
 \end{center}
\bigskip
 The maun objective of this paper is to study the case of flag supermanifolds of maximal type corresponding to the Lie superalgebra $\mathfrak{osp}_{2m+1|2n}(\mathbb{C})$.

   The Lie superalgebra $\mathfrak{osp}_{m|2n}(\mathbb C)$ is the linear Lie superalgebra that annihilates a non-degenerate  even symmetric bilinear form in $\mathbb C^{m|2n}$. (See \cite{Kac} for details about Lie superalgebras.) Let us take two $r$-tuples $k=(k_0,\ldots,k_r)$ and $l=(l_0,\ldots,l_r)$ of non-negative integers such that
   \begin{equation}\label{eq conditions on k_i and l_j}
   \begin{split}
   0\le k_r\le\ldots\le k_0= m&,\quad 0\le l_r\ldots\le l_0= n,\\
   \quad 0 <
   k_r +l_r <\ldots& < k_0 + l_0 = m+n.
   \end{split}
   \end{equation} 
   We denote by $\mathbf
 F_{k|l}$ the flag supermanifold of type $k|l$ corresponding to the Lie superalgebra $\mathfrak{gl}_{m|n}(\mathbb C)$. This is a super-analog of usual flag manifold $\mathbf
 F_{k}$ of type $k=(k_0,\ldots,k_r)$ in $\mathbb C^{m}$.   Further we denote  by $\mathbf {IF}_{k|l}$  the isotropic flag supermanifold in $\mathbb C^{m|2n}$ corresponding to the Lie superalgebra $\mathfrak{osp}_{m|2n}(\mathbb C)$. We set $k'=(k_1,\ldots,k_r)$ and $l'= (l_1,\ldots,l_r)$. It was proven in \cite{ViGL in JA} that under some restriction on the flag type we have $\mathfrak{v}(\mathbf{F}_{k'|l'})\simeq
 \mathfrak{pgl}_{k_1|l_1}(\mathbb C)$. Here we put  $\mathfrak{pgl}_{m|n}(\mathbb C):=\mathfrak{gl}_{m|n}(\mathbb C)/ \mathfrak
 z(\mathfrak{gl}_{m|n}(\mathbb C))$, where $\mathfrak
 z(\mathfrak{gl}_{m|n}(\mathbb C))$ is the center of $\mathfrak{gl}_{m|n}(\mathbb C)$.  The main result of our paper is the following.

\medskip

\t\label{teor main osp} {\sl  Let $r>1$, $k_1\geq 1$ and
	$l_1\geq 1$, the conditions (\ref{eq_condition no functions}) hold, $\mathfrak{v}(\mathbf{F}_{k'|l'})\simeq \mathfrak {pgl}_{k_1|l_1}(\mathbb
	C)$ and $k'\ne (k_1-1,\ldots, k_1-1,0,\ldots,0)$. Then 
	$\mathfrak{v}(\mathbf F_{k|l}(\mathfrak{osp}_{2k_1-1|2l_1}(\mathbb
	C)))\simeq \mathfrak{osp}_{2k_1-1|2l_1}(\mathbb C)$.} $\Box$

\medskip

\bigskip

\textbf{Acknowledgements:} E.~V. was partially  supported by SFB TR 191 and by the Universidade Federal de Minas Gerais.

\section{Lie superalgebra $\mathfrak {osp}_{2m+1|2n}(\mathbb C)$}

More information about  Lie supergroups  and Lie superalgebras can be found for instance \cite{Kac,ViLieSupergroup}. 
 The Lie superalgebra $ \mathfrak {osp}_{2m+1|2n}(\mathbb C)$ is the Lie sub-superalgebra in $\mathfrak{gl}_{2m+1|2n}(\mathbb C)$ that annihilates a non-degenerate even symmetric bilinear form $\beta$ in $\mathbb
C^{2m+1|2n}$. If the matrix $\Gamma$ of $\beta$ in the standard basis of  $\mathbb C^{2m+1|2n}$ is given by
\begin{equation}
\label{matr_grama_osp}
 \Gamma=\left(
\begin{array}{ccccc}
0&E_m & 0 & 0 & 0 \\
E_m& 0 & 0 & 0 & 0 \\
0&0 & 1 & 0 & 0 \\
0&0 & 0 & 0 & E_n \\
0&0 & 0 & -E_n & 0 \\
\end{array}
\right),
\end{equation}
then the Lie superalgebra $\mathfrak {osp}_{2m+1|2n}(\mathbb C)$ is defined by the equation $M^{ST}\Gamma+ \Gamma M =0$, 
where
$$
\left(
\begin{array}{cc}
M_{11} & M_{12} \\
M_{21} & M_{22} \\
\end{array}
\right)^{ST}=\left(
\begin{array}{cc}
M_{11}^T & M_{21}^T \\
-M_{12}^T & M_{22}^T \\
\end{array}
\right).
$$
and $T$ is the usual transposition. Explicitly we have 
\begin{equation}
\label{matr_osp_nechet} \mathfrak {osp}_{2m+1|2n}(\mathbb C)=
\left\{ \left(
\begin{array}{ccccc}
A_{11}&A_{12}&G_1&C_{11}& C_{12}\\
A_{21} &-A_{11}^T & G_{2} & C_{21} & C_{22} \\
-G_2^T &   -G_1^T & 0 & G_3 & G_4 \\
-C_{22}^T & -C_{12}^T & -G^T_{4} & B_{11} & B_{12} \\
C_{21}^T &   C_{11}^T & G_{3}^T & B_{21} & -B_{11}^T \\
\end{array}
\right),
\begin{array}{l}
A_{21}^T=-A_{21},\\A_{12}^T=-A_{12},\\
B_{12}^T=B_{12}, \\ B_{21}^T=B_{21}
\end{array}
\right\}.
\end{equation}
Here $A_{11}$, $B_{11}$ are square complex matrices of size $m$ and $n$, respectively, $G_1$ and $G_2$ are complex matrices of size $m\times 1$. The center $\mathfrak z(\mathfrak {osp}_{2m+1|2n}(\mathbb
C))$ of $\mathfrak {osp}_{2m+1|2n}(\mathbb
C)$ is trivial. The corresponding Lie supergroup  we will denote by $\operatorname{OSp}_{2m+1|2n}(\mathbb C)$. This is a sub-supermanifold in $\operatorname{GL}_{2m+1|2n}(\mathbb C)$ that is given by the following equation:
$$
\left(
\begin{array}{cc}
X_{11} & X_{12} \\
X_{21} & X_{22} \\
\end{array}
\right)^{ST}\Gamma \left(
\begin{array}{cc}
X_{11} & X_{12} \\
X_{21} & X_{22} \\
\end{array}
\right)=\Gamma.
$$

\section{Isotropic flag supermanifolds}

 More information about the theory of (complex) supermanifolds can be found in \cite{BL,Fioresi,Leites,Man}.  We will denote a complex-analytic supermanifold by $\mathcal{M} = (\mathcal{M}_0,{\mathcal O}_{\mathcal{M}})$, where $\mathcal{M}_0$ is the underlying complex-analytic manifold and ${\mathcal O}_{\mathcal{M}}$ is the structure sheaf of $\mathcal{M}$. 
 The Lie superalgebra of global sections of the sheaf $\mathcal Der\,(\mathcal O_{\mathcal M})$ is called the {\it Lie superalgebra of holomorphic vector fields} on $\mathcal M$. 
  We will denote this Lie superalgebra by $\mathfrak v(\mathcal M)$. The Lie superalgebra $\mathfrak v(\mathcal M)$ is finite dimensional if $\mathcal M_0$ is compact. A proof of this statement can be found for instance in \cite{ViGL in JA}.

\medskip

\noindent{\bf $\mathfrak {gl}_{m|n}(\mathbb C)$-flag supermanifolds.} Let us recall briefly a construction of an atlas on the supermanifold  $\mathbf F_{k|l}$, where $k=(k_0,\ldots,k_r)$ and $l=(l_0,\ldots,l_r)$ such that (\ref{eq conditions on k_i and l_j}) holds. More details can be found in \cite{ViGL in JA,Vi osp and pi}. Recall that the underlying space of the supermanifold $\mathbf F_{k|l}$ is the product  $\mathbf
F_{k}\times \mathbf F_{l}$ of two classical flag manifolds of types $k$ and $l$, respectively. Consider two subsets $I_{s\bar
0}\subset\{1,\ldots,k_{s-1}\}$  and $I_{s\bar 1}\subset\{1,\ldots,l_{s-1}\}$
  such that $|I_{s\bar 0}| = k_s,$ and $|I_{s\bar 1}| = l_s$,  where $s = 1,\ldots,r$. We set $I_s:=(I_{s\bar 0},I_{s\bar 1})$ and $I := (I_1,\ldots,I_r)$. To any set $I_s$  we assign the following $(k_{s-1} + l_{s-1})\times (k_s + l_s)$-matrix 
\begin{equation}\label{eq local chart general}
Z_{I_s}
=
\left(
\begin{array}{cc}
X_s & \Xi_s\\
\H_s & Y_s \end{array} \right), \ \ s=1,\dots,r.
\end{equation}
The matrices $X_s=(x^s_{ij})$ and $Y_s=(y^s_{ij})$ in (\ref{eq local chart general}) are of the size $(k_{s-1}\times k_s)$ and $(l_{s-1}\times l_s)$, respectively. Moreover, we assume that the matrix $Z_{I_s}$ contains the identity submatrix $E_{k_s+l_s}$ of size $(k_s+l_s)\times (k_s+l_s)$ in the lines with numbers $i\in I_{s\bar 0}$ and $k_{s-1} + i,\; i\in I_{s\bar 1}$. 
The matrices $(Z_{I_s})$, where $s=1,\dots,r$,  define the superdomain $\mathcal U_I$ with even coordinates $x^s_{ij}$ and $y^s_{ij}$, and odd coordinates $\xi^s_{ij}$ and $\eta^s_{ij}$, where $\xi^s_{ij}$ and $\eta^s_{ij}$ are non-trival entries of $\Xi_s$ and $\H_s$, respectively. In the case $r = 1$ the supermanifold $\mathbf F_{k|l}$ is called the {\it super-Grass\-mannian}. Sometimes the notation $\mathbf {Gr}_{m|n,k|l}$ is used for $\mathbf F_{k|l}$ with $r=1$ in the literature. 

The supermanifold $\mathbf F_{k|l}$ is homogeneous, since it  possesses a transitive action of the Lie supergroup $\GL_{m|n}(\mathbb C)$, see \cite{ViGL in JA,Viholom} for details and definitions. 
 Let us recall the  definition of this action  in our atlas. Let 
 $$
 L= 
 \left(
 \begin{array}{cc}
 L_{11} & L_{12} \\
 L_{21} & L_{22} \\
 \end{array}
 \right)
 $$
 be a coordinate matrix of the Lie supergroup $\GL_{m|n}(\mathbb C)$. Then the action of  $\GL_{m|n}(\mathbb C)$ on 
 $\mathbf {F}_{k|l}$ in our coordinates is given by:
\begin{equation}\label{eq action of Q}
\begin{aligned} 
&(L,(Z_{I_1},\ldots,Z_{I_r})) \longmapsto
(\tilde Z_{J_1},\ldots,\tilde Z_{J_r}), 
\end{aligned} 
\end{equation}
where $\tilde Z_{J_1} =
LZ_{I_1}C_1^{-1}$ and $\tilde Z_{J_s} = C_{s-1}Z_{I_s}C_s^{-1}$. 
Here $C_1$ is the invertible submatrix in $LZ_{I_1}$ that consists of the lines with numbers $i\in J_{1\bar 0}$ and $k_0+i$, where  $i\in J_{1\bar 1}$, and $C_s,$ where $s\ge 2$, is the invertible submatrix in $C_{s-1}Z_{I_s}$ that consists of the lines with numbers $i\in J_{s\bar 0}$ and $k_{s-1}+i$, where
$i\in J_{\bar 1s}$. This Lie supergroup action induces the Lie superalgebra homomorphism  $\nu:\mathfrak {gl}_{m|n}(\mathbb C)\to\mathfrak v(\mathbf {F}_{k|l}).
$

\medskip

\noindent{\bf  $\mathfrak
	{osp}_{m|2n}(\mathbb C)$-flag supermanifolds.}
 The supermanifold $\mathbf {IF}_{k|l}$, where $l_0$ is even, is a sub-supermanifold of  $\mathbf F_{k|l}$ that is given in local coordinates $(\ref{eq local chart general})$ by the following equation:
\begin{equation}
\begin{split}
\left(
\begin{array}{cc}
X_1 & \Xi_1\\
\H_1 & Y_1 \end{array} \right)^{ST} \Gamma \left(
\begin{array}{cc}
X_1 & \Xi_1\\
\H_1 & Y_1 \end{array} \right)=0,
\end{split}
\end{equation}
where $\Gamma$ is defined by (\ref{matr_grama_osp}) and $ST$ is the super-transposition as above. There are transitive actions $\mu$ of the Lie supergroup  $\OSp_{m|2n}(\mathbb C)$ on $\mathbf {IF}_{k|l}$. It is given by Formulas  $(\ref{eq action of Q})$, if we replace  $L$ by a coordinate matrix of the Lie supergroup $\OSp_{m|2n}(\mathbb C)$. This action induces the Lie superalgebras homomorphism
$$
\mu:\mathfrak
{osp}_{m|2n}(\mathbb C)\to\mathfrak v(\mathbf {IF}_{k|l}).
$$
In the case 

{\bf (1)} $m=2k_1+1$ and $l_0=2l_1$ or 

{\bf (2)} $m=2k_1$ and $l_0=2l_1$, 

\noindent we say that the supermanifold $\mathbf {IF}_{k|l}$ has {\it maximal type}. In this paper we study the case $m=2k_1+1$ and $l_0=2l_1$. The case {\bf (2)} can be found in \cite{Vi osp and pi}.

\section{Vector fields on a superbundle}

Recall that a morphism of supermanifolds is a morphism of the corresponding ringed spaces that preserves the $\mathbb Z_2$-grading of the structure sheaves. Let $\pi=(\pi_0,\pi^*): \mathcal M\to \mathcal N$ be a morphism of supermanifolds. As in the classical situation, a vector field $v\in\mathfrak
v(\mathcal M)$ is called {\it projectible} with respect to $\pi$, if there exists a vector field $v_1\in\mathfrak v(\mathcal N)$ such that $\pi^*(v_1(f))=v(\pi^*(f))$ for all $f\in \mathcal O_{\mathcal N}$. In the case when 
$\pi=(\pi_0,\pi^*): \mathcal M\to \mathcal N$ is the projection of a superbundle, the vector field $v_1$ is unique if it exists.  A vector field $v\in\mathfrak v(\mathcal M)$ is called {\it vertical}, if it is projectible and its projection $v_1$ is equal to $0$. If $\mathcal S$ is a supermanifold, then the global sections $\mathcal O_{\mathcal S}(\mathcal S_0)$ of the structure sheaf $\mathcal O_{\mathcal S}$ are called {\it holomorphic functions} on $\mathcal S$.  It was shown in \cite{Bash} that if  
$\pi: \mathcal M 	\to \mathcal B$ is the projection of a superbundle with fiber $\mathcal S = (\mathcal S_0, \mathcal O_{\mathcal S})$ and $\mathcal O_{\mathcal S}(\mathcal S_0) = \mathbb C$, then any global holomorphic vector field on $\mathcal M$ is projectible. In other words there exists a homomorphism of Lie superalgebras 
$$
\mathcal  P: \mathfrak v(\mathcal M) \to \mathfrak v(\mathcal B).
$$

For computation of the Lie superalgebra of holomorphic vector fields on an isotropic flag supermanifold we will use the following fact. For $r> 1$ the isotropic flag supermanifold $\mathbf {IF}_{k|l}$ is a superbundle  with the base space that is isomorphic to the isotropic super-Grassmannian $\mathbf {IF}_{k_0,k_1|l_0,l_1}$  and with the fiber that is isomorphic to $\mathbf F_{k'|l'}$.  In our local coordinates the bundle projection that is denoted by $\pi$ is given by $(Z_1,Z_2,\ldots Z_n) \longmapsto (Z_1).
$ Moreover, from Formulas (\ref{eq action of Q}) we can deduce that the projection $\pi$ is equivariant with respect to the action of the Lie supergroup $\OSp_{m|2n}(\mathbb C)$  on $\mathbf {IF}_{k|l}$. Further, the following holds.

\medskip
\t\cite{Viholom}\label{teor constant functions} {\sl Consider the flag supermanifold $\mathcal M=\mathbf{F}_{k|l}$.  Assume that 
	\begin{equation}\label{eq_condition no functions}
	\begin{split}
	(k|l)&\ne 
	(m,\ldots, m, k_{s+2},\ldots, k_r| l_1,\ldots, l_s,0,\ldots,0),\\
	(k|l)&\ne (k_1,\ldots, k_s,0,\ldots,0 |n,\ldots, n, l_{s+2},\ldots,
	l_r),
	\end{split}
	\end{equation}
	for any $s\geq 0$. Then $\mathcal O_{\mathcal M}(\mathcal M_0) = 	\mathbb{C}.$ 	Otherwise
	$\mathcal O_{\mathcal M}(\mathcal M_0) = \bigwedge(mn)$, where $\bigwedge(mn)$ is the Grassmann algebra with $mn$ generators.$\Box$ 
}

\medskip

Assume that the fiber $\mathbf F_{k'|l'}$ satisfies the conditions of Theorem \ref{teor constant functions}. Then all holomorphic vector fields on  $\mathbf {IF}_{k|l}$ are projectible and we have the following Lie algebra homomorphism 
\begin{equation}\label{eq mathcal P}
\mathcal P: \mathfrak v(\mathbf {IF}_{k|l}) \to \mathfrak v(\mathbf {IF}_{k_0,k_1|l_0,l_1}).
\end{equation}

As it was noticed above,  the action  (\ref{eq action of Q}) of the Lie supergroup $\operatorname{OSp}_{m|2n}(\mathbb C)$ induces the following Lie superalgebras homomorphism:
$$
\mu: \mathfrak {osp}_{m|2n}(\mathbb C)\to
\mathfrak{v}(\mathbf {IF}_{k|l}).
$$
The Lie superalgebra $\mathfrak{v}(\mathbf {IF}^o_{m,k_1|n,l_1})$ of holomorphic vector fields on a connected component $\mathbf {IF}^o_{m,k_1|n,l_1}$ of the isotropic  super-Grass\-mannian $\mathbf {IF}_{m,k_1|n,l_1}$ of maximal type was calculated in  \cite{oniosp}. 

\medskip
\t[Onishchik-Serov]\label{teor super=Grassmannians} {\sl Let $r=1$.	
	Assume that $m=2k_1$ and $n=2l_1$, i.e the super-Grassmannian $\mathbf {IF}_{m,k_1|n,l_1}$ is of maximal type. If $k_1\geq 1$ and $l_1\geq 1$, then the homomorphism 
	$$
	\mu: \mathfrak {osp}_{m|n}(\mathbb C)\to
	\mathfrak{v}(\mathbf {IF}^o_{m,k_1|n,l_1})
	$$
	is an isomorphism.

}

\medskip

If $k_1\geq
1$ and $l_1\geq 1$, the super-Grassmannian $\mathbf {IF}_{m-1,k_1-1|n,l_1}$  is isomorphic to $\mathbf {IF}^o_{m,k_1|n,l_1}$, this is to a connected component of the super-Grassmannian $\mathbf {IF}_{m,k_1|n,l_1}$. A proof of this result can be found below. Therefore the Lie superalgebra of vector fields on  the isotropic  super-Grass\-mannian $\mathbf {IF}_{m-1,k_1-1|n,l_1}$ is also isomorphic to $\mathfrak {osp}_{m|n}(\mathbb C)$.

\section{The Borel-Weil-Bott Theorem}

To compute the Lie superalgebra of vector fields on isotropic flag supermanifolds we will use the Borel-Weil-Bott Theorem.  Details can be found for example in \cite{ADima}.
  Consider the isotropic super-Grassmanian of maximal type $\mathbf {IF}_{2s+1,s|2n,n}$. The underlying manifold of $\mathbf {IF}_{2s+1,s|2n,n}$ is a homogeneous space $G/P$, where $G\simeq \operatorname{SO}_{2s+1}(\mathbb C)\times \operatorname{Sp}_{2n}(\mathbb
 C)$ is the underlying space of the Lie supergroup 
$\operatorname{OSp}^o_{2s+1|2n}(\mathbb C)$, where $\operatorname{OSp}^o_{2s+1|2n}(\mathbb C)$ is the connected component of $\operatorname{OSp}_{2s+1|2n}(\mathbb C)$ that contains identity element, and $P$ is the parabolic subgroup in $G$ that contains all matrices in the following form:
\begin{equation}
\label{P_osp_m-1}
\left(
\begin{array}{ccccc}
A_1 & 0& 0& 0& 0\\
C_1& (A_1^T)^{-1}& G& 0& 0\\
H&0&1& 0& 0\\
0&0&0& A_2 & 0\\
0&0&0& C_2& (A_2^T)^{-1}
\end{array}
\right).
\end{equation}
Here 
$A_1\in \GL_{s}(\mathbb C)$ and $A_2\in \GL_{n}(\mathbb C)$.
The reductive part $R$ of $P$ has the form
$R\simeq
\GL_{s}(\mathbb C)\times \GL_{n}(\mathbb C)$.

Further, denote by $\mathfrak {osp}_{2s+1|2n}(\mathbb C)_{\bar 0} =\mathfrak{so}_{2s+1}(\mathbb
C)\oplus\mathfrak {sp}_{2n}(\mathbb C)$ the even part of the Lie superalgebra $\mathfrak {osp}_{2s+1|2n}(\mathbb C)$. In the Lie algebra $\mathfrak {osp}_{2s+1|2n}(\mathbb C)_{\bar 0}$ we fix the following data: a Cartan subalgebra $\mathfrak t=\mathfrak
t_1 \oplus \mathfrak t_2$, a system of positive roots $\Delta^+ = \Delta^+_1 \cup \Delta^+_2$ and a system of simple roots $\Phi= \Phi_1\cup \Phi_2$. More precisely, we have for $s\geq 1$ and $n\geq 1$
 \begin{align*}
\mathfrak t_1&=
\{\operatorname{diag}(a_1,\dots,a_{s}, -a_1,\dots,-a_{s}),0\},\quad
\mathfrak t_2=
\{\operatorname{diag}(b_1,\dots,b_n,-b_1,
\dots,-b_n)\};\\
 \Delta^+_1&=\{\mu_i-\mu_j,\,\mu_i+\mu_j, \,\,i<j, \,\, \mu_i,\},
\quad
\Delta^+_2=\{\lambda_p-\lambda_q,\,\,
\,\,p<q, \,\, \lambda_p+\lambda_q,\,\, \,\,p\leq q \};\\
\Phi_1&= \{\alpha_1,...,
\alpha_{s}\}, \,\,\, \alpha_i=\mu_i-\mu_{i+1},
\,\,i=1,\ldots, s-1, \,\,\alpha_{s}=\mu_{s};\\
\Phi_2&= \{\beta_1,..., \beta _{n}\},
\,\,\, \beta_j=\lambda_j-\lambda_{j+1},\,\,j=1,\ldots, n-1, \,\,
\beta_n=2\lambda_n.
\end{align*}
Note that in the case $s=1$ we have an isomorphism $\mathfrak{so}_{3}(\mathbb C)\simeq \mathfrak{sl}_{2}(\mathbb
C)$. Similarly in the case  $n=1$, we have $\mathfrak {sp}_{2}(\mathbb C)\simeq \mathfrak{sl}_{2}(\mathbb C)$. However in both cases we also can use the root system $\Delta^+$.

 Denote by $\mathfrak t^*(\mathbb R)$
the real subspace in $\mathfrak t^*$ 
spanned by $\mu_j$ and $\lambda_i$ equipped with  the scalar product $( \,,\, )$ such that the vectors  $\mu_j,\lambda_i$ form an orthonormal basis. An element $\gamma\in \mathfrak t^*(\mathbb R)$ is called {\it dominant} if $(\gamma, \alpha)\ge 0$ for all positive roots $\alpha \in \Delta^+$.

Let $\psi$ be a representation of the group $P$ in a finite dimensional complex vector space $E$. To such a representation we can assign a homogeneous holomorphic vector bundle $\mathbf E_{\psi}\to G/P$ with the fiber $E=(\mathbf E_{\psi})_{P}$ at the point $P$.
 We denote by $\mathcal E_{\psi}$ the shaef of holomorphic sections of $\mathbf E_{\psi}$.

\medskip

\t [Borel-Weil-Bott]. \label{teor borel} {\sl Assume that the representation	$\psi: P\to \GL(E)$ is completely reducible and $\rho_1,..., \rho_s$
	are highest weights of $\psi|R$. Then the $G$-module $H^0(G/P,\mathcal E_{\psi})$ is isomorphic to the sum of irreducible $G$-modules with highest weights $\rho_{i_1},..., \rho_{i_t}$, where 
	$\rho_{i_a}$ are dominant highest weights.

}

\medskip

We will use Theorem \ref{teor borel} for the following holomorphic vector bundle. Consider a superbundle $\pi: \mathcal M\to \mathcal B$ with fiber $\mathcal S$. Denote by $\mathcal W$ the following sheaf on $\mathcal B_0$: to any open set $U\subset \mathcal B_0$ we assign
the Lie superalgebra of all vertical vector fields on the supermanifold $(\pi_0^{-1}(U),\mathcal
O_{\mathcal M})$. If the base space $\mathcal S_0$ of the fiber $\mathcal S$ is compact, then $\mathfrak v(\mathcal S)$ is finite dimensional  and $\mathcal W$ is a locally free sheaf of $\mathcal O_{\mathcal B}$-modules with
$\dim\mathcal W= \dim\mathfrak v(\mathcal S)$, see \cite{ViGL in JA} for details. In \cite{ViGL in JA} a description of the corresponding to $\mathcal W$ graded sheaf $\widetilde {\mathcal W}$ was obtained. It is defined in the following way:
\begin{equation}\label{eq tilda W def}
	\widetilde {\mathcal W}=\bigoplus_{p\ge 0}\widetilde{\mathcal W}_{p}, \quad
	\text{where}\quad
	\widetilde{\mathcal W}_{p}=\mathcal W_{(p)}/\mathcal
	W_{(p+1)}.
\end{equation}
Here $\mathcal W_{(p)}=
\mathcal J^p\mathcal W$ and $\mathcal J$ is the sheaf of ideals in ${\mathcal O}_{\mathcal B}$ generated by odd elements. 
Clearly, $\widetilde {\mathcal W}$ is a $\mathbb Z$-graded 
sheaf of $\mathcal F_{\mathcal B_0}$-modules, where  $\mathcal F_{\mathcal B_0}$ is the structure sheaf of the underlying manifold $\mathcal B_0$. If the base space $\mathcal S_0$ is compact, then $\widetilde{\mathcal W}_0$ is a locally free sheaf of $\mathcal F_{\mathcal B_0}$-modules and any fiber of the corresponding vector bundle is isomorphic to $\mathfrak v(\mathcal S)$, see \cite{ViGL in JA} for details.

\section{Vector fields on $\mathbf
	{IF}_{k|l}$, case $m=2k_1-1$}

In this section we will compute the Lie superalgebra of holomorphic vector fields on $\mathbf {IF}_{k|l}$ of type
$$
k= (2k_1-1,k_1-1,k_2,\ldots,k_r),\quad l= (2l_1,l_1,l_2,\ldots,l_r),
$$
where $r>1$. We also put $m=2k_1-1$ and $n=2l_1$.  In other words, we assume that $m$ is odd and that $\mathbf
{IF}_{k|l}$ is an isotropic flag supermanifold of maximal type.  As we have seen above the isotropic flag supermanifold $\mathcal M:=\mathbf {IF}_{k|l}$ is a superbundle. We denote by $\mathcal B = \mathbf {IF}_{2k_1-1,k_1-1|2l_1,l_1}$ and by $\mathcal S =\mathbf{F}_{k'|l'}$ its base space and its fiber, respectively. Here 
$$
k'=(k_1-1,k_2,\ldots,k_r),\quad l'= (l_1,l_2,\ldots,l_r).
$$ 
Note that $\mathcal S$ is a usual flag supermanifold.

In the previous sections we considered the homomorphism $\mu:
\mathfrak{osp}_{m|n}(\mathbb C)\to \mathfrak v(\mathcal M)$.  Denote
$\mu_{\mathcal B}: \mathfrak{osp}_{m|n}(\mathbb C)\to\mathfrak v(\mathcal B)$ the corresponding homomorphism for the super-Grassmannian. Since the projection of the superbundle $\mathcal M$ is equivariant, these homomorphisms satisfy the relation $\mu_{\mathcal B} =
\mathcal{P}\circ\mu$, where $\mathcal{P}$ is defined by (\ref{eq mathcal P}). From Theorem \ref{teor super=Grassmannians} it follows that $\mu_{\mathcal B}$ is not surjective. In the next sections we will show that  $\Ker \mathcal{P}=\{0\}$ and $\Im \mathcal{P} = \mathfrak{osp}_{m|n}(\mathbb C)$ and we will conclude that $\mu:
\mathfrak{osp}_{m|n}(\mathbb C)\to \mathfrak v(\mathcal M)$ is an isomorphism.

\subsection{Computation of $\Ker \mathcal{P}$}

On the base space $\mathcal B$  of the superbundle $\mathbf
{IF}_{k|l}$ in the case $m=2k_1-1$ we can choose the following chart:
\begin{equation}
	\label{okrestnost_osp_2m-1} 
	Z_{I_1}=
	\left(
	\begin{array}{cc}
		Z_1 & \mathcal{Z}_1 \\
		E_{k_1-1} & 0 \\
		X_1 & \Xi_1 \\
		\H_1 & Y_1 \\
		0 & E_{l_1} \\
	\end{array}
	\right),
\end{equation}
where $Z_1=(z^1_{ij})$, $X_1=(x^1_{i})$, $Y_1=(y^1_{ij})$ are matrices with even coordinates of sizes 
$(k_1-1)\times (k_1-1)$, $1\times
(k_1-1)$ and $l_1\times l_1$, respectively,
$\mathcal{Z}_1=(\zeta^1_{ij})$, $\Xi_1=(\xi^1_{i})$,
$\H_1=(\eta^1_{ij})$ are matrices with odd coordinates. As above $E_p$ is the identity matrix of size $p$. Moreover the following relations hold:
$$
\begin{array}{l}
Z_1+Z_1^T+X_1^TX_1=0,\quad
\mathcal{Z}_1^T+\Xi_1^TX_1+\H_1=0, \quad
\Xi_1^T\Xi_1+Y_1-Y_1^T=0.
\end{array}
$$
We can choose elements from $X_1$,
$\Xi_1$, $\H_1$ and elements from $Z_1$ and $Y_1$ strictly under the diagonal and  under the diagonal, respectively,  as independent coordinates. We choose some set of indexes $I_t$ and the coordinate matrices $Z_{I_t}$, where $t>1$, to complete the chart $Z_{I_1}$ to a chart on $\mathbf {IF}_{k|l}$. The obtained chart we denote by $\mathcal U_I$. We need the following lemma.

\medskip

\lem \label{lem fields osp_2m-1} {\it The following vector fields  on $\mathbf {IF}_{k|l}$ written in coordinates of $\mathcal U_I$ 
	$$
	\frac{\partial}{\partial \eta^1_{ab}},\quad
	h_i=\frac{\partial}{\partial \xi^1_{i}}-
	\sum_jx^1_j\frac{\partial}{\partial \eta^1_{ij}}- \sum_{j\leq
		i}\xi^1_j\frac{\partial}{\partial y^1_{ij}}
	$$
are fundamental. In other words these vector fields are induced by the natural action $\mu:
\mathfrak{osp}_{m|n}(\mathbb C)\to \mathfrak v(\mathcal M)$.

}
\medskip

\noindent {\it Proof}. The vector field $h_i$, where $i=1,\ldots,l_1$, corresponds to the Lie sub-supergroup 
$$
\gamma(\tau)=\exp(-\tau
E_{2k_1-1+i,2k_1-1}+ \tau E_{2k_1-1,2k_1-1+l_1+i}),
$$
where $\tau$ is an odd parameter and $E_{ij}$ is the matrix with $1$ in the position $(ij)$ and $0$ otherwise. Indeed, $\gamma(\tau)$ transforms the coordinate matrices in the following way
$$
\left(
\begin{array}{cc}
Z_1 & \mathcal{Z}_1 \\
E_{k_1-1} & 0 \\
X_1 & \Xi_1 \\
\H_1 & Y_1 \\
0 & E_{l_1} \\
\end{array}
\right)\mapsto \left(
\begin{array}{cc}
Z_1 & \mathcal{Z}_1 \\
E_{k_1-1} & 0 \\
X_1 & \tilde{\Xi}_1 \\
\tilde{\H}_1 & \tilde Y_1 \\
0 & E_{l_1} \\
\end{array}
\right), \quad
 Z_{I_s}\mapsto Z_{I_s},\,\, s\geq 2,
$$
where
$$
\tilde{\Xi}_1=\Xi_1+\tau F_i, \quad\tilde{\H}_1 =\H_1-\tau
F_i^TX_1,\quad \tilde Y_1 =Y_1-\tau F_i^T\Xi_1.
$$
Here $F_i=(0,\ldots,1,\ldots,0)$ has $1$ on the $i$-th place. The result follows. The case of the vector fields $\frac{\partial}{\partial \eta^1_{ij}}$ is similar.$\Box$

\medskip

Recall that we have the following exact sequence of sheaves.
$$
0\to \mathcal{W}_{(1)} \longrightarrow \mathcal{W}_{(0)} \longrightarrow \widetilde{\mathcal{W}}_{0} \to 0. 
$$
Hence we have an exact sequence of $0$-cohomology groups
$$
0\to H^0(\mathcal B_0, \mathcal{W}_{(1)}) \longrightarrow H^0(\mathcal B_0,\mathcal{W}_{(0)}) \longrightarrow H^0(\mathcal B_0,\widetilde{\mathcal{W}}_{0}) . 
$$
We will need the following lemma. 

\medskip

\lem \label{lem kerosp_2m-1} {\it Assume that  $\Ker \mathcal{P}\ne \{0\}$. Then we have
	$$
	H^0(\mathcal B_0, \mathcal{W}_{(1)}) \subsetneq H^0(\mathcal B_0,\mathcal{W}_{(0)}).
	$$
	
}

\medskip

\noindent{\it Proof.} Let us take $v\in\Ker \mathcal{P}\setminus
\{0\}$. Recall that $\mathcal{W}$ is a locally free sheaf of $\mathcal O_{\mathcal B}$-modules and
$\dim\mathcal W= \dim\mathfrak v(\mathcal S)<\infty$. Let us take a basis $(v_j)$ in $\mathfrak v(\mathcal S)$. Then  in our chart $\mathcal U_I$ we can write the  vertical vector field $v$ in the form  $v=\sum f_jv_j$, where $f_j\in
\mathcal{O}_B$. Clearly, $\Ker \mathcal{P}$ is an ideal in $\mathfrak v(\mathcal M)$. Note that we may assume that $f_j$,  are independent on coordinates $\eta_{ab}^1$. 
 Indeed, if some $f_s$ is dependent on $\eta_{ab}^1$, then the commutator $[\frac{\partial}{\partial \eta^1_{ab}}, v]\in\Ker \mathcal{P}$ is not $0$ and does not contain $\eta_{ab}^1$ anymore.

Further,  let us write $f_j$ in the form:
$$
f_j=f^j_{q(j)}+f^j_{q(j)+1}+\ldots,
$$
where $f^j_{q(j)+k}$ is a homogeneous polynomial of degree $q(j)+k$ in  $\xi_i^1$. We put $q=\min_j\{q(j)\}$. Then
$$
v\in H^0(\mathcal B_0,\mathcal{W}_{(q)}) \backslash
H^0(\mathcal B_0,\mathcal{W}_{(q+1)}).
$$
If $q=0$, everything is proven. Assume that $q>0$ and that  $v_q:=\sum\limits_{q(j)=q}f^j_{q(j)}v_j$ depends on a coordinate $\xi_i^1$ for some $i$. Denote $v_1:= v-v_q$. We have
\begin{align*}
[h_i,v]= [h_i,v_q]+[h_i,v_1], \quad [h_i,v_q]\in \mathcal{W}_{(q-1)}(\mathcal B_0)\quad \text{and} \quad [h_i,v_1]\in \mathcal{W}_{(q)}(\mathcal B_0). 
\end{align*}
Hence,
$$
[h_i, v]\in H^0(\mathcal B_0,\mathcal{W}_{(q-1)}) \backslash
H^0(\mathcal B_0,\mathcal{W}_{(q)}).
$$ 
This completes the proof. $\Box$

\medskip

Denote by $\widetilde{\mathcal W}_0$ the locally free sheaf of $\mathcal F_{\mathcal B_0}$-modules on $ \mathcal B_0=G/P$, see Section $5$. 
Let us compute the representation $\psi$ of
$P$ in the fiber $W$ of  $\widetilde{\mathcal W}_0$  over the point 
$o:=P$.
Let us compute first the action of $P$ on the fiber  $\mathcal S_o=\mathcal S$. The point $o$ is given in the chart $\mathcal U_I$ by the equations
$X_1=0$, $Y_1=0$, $Z_1=0$, $\Xi_1=0$, $\H_1=0$, $\mathcal{Z}_1=0$.
The action of $P$ at $o$  in coordinates $Z_{I_1}$ is given by
$$
\left(
\begin{array}{ccccc}
A_{1} & 0              & 0 & 0 & 0 \\
C_{1} & (A_{1}^T)^{-1} & G & 0 & 0 \\
H & 0 & 1 & 0 & 0 \\
0 & 0 & 0 & A_{2} & 0 \\
0 & 0 & 0 & C_{2} & (A_{2}^T)^{-1} \\
\end{array}
\right)
\left(
\begin{array}{cc}
0 & 0 \\
E_{k_1-1} & 0 \\
0 & 0 \\
0 & 0 \\
0 & E_{l_1} \\
\end{array}
\right)= \left(
\begin{array}{cc}
0 & 0 \\
(A_{1}^T)^{-1} & 0 \\
0 & 0 \\
0 & 0 \\
0 & (A_{2}^T)^{-1} \\
\end{array}
\right).
$$
Hence, in coordinates $Z_{I_2}$ we have
$$
\left(
\begin{array}{cc} (A_{1}^T)^{-1} & 0\\ 0 &  (A_{2}^T)^{-1}\end{array}
\right) \left(
\begin{array}{cc}  X_2 & \Xi_2 \\ \H_2 & Y_2\end{array} \right) =
\left(
\begin{array}{cc}  (A_{1}^T)^{-1}X_2  &
(A_{1}^T)^{-1}\Xi_2 \\
(A_{2}^T)^{-1} \H_2 &
(A_{2}^T)^{-1}Y_2
\end{array} \right).
$$
Other coordinate matrices  $Z_{I_p}$, $p>2$, are transformed accordingly.   Therefore, the action of $P$ on $\mathcal S_o$ coincides with the composition of the standard action of $\GL_{k_1-1}(\mathbb C)\times \GL_{l_1}(\mathbb C)$ and $(A_1,A_2) \to ((A_{1}^T)^{-1}, (A_{2}^T)^{-1})$.

Assume that
$$
\mathfrak v(\mathcal S)\simeq\mathfrak{pgl}_{k_1-1|l_1}(\mathbb
C)=\left\lbrace\left(
\begin{array}{cc}
Z_1 & Q_1  \\
Q_2 & Z_2  \\
\end{array}
\right)+<E_{k_1+l_1-1}>\right\rbrace,
$$
where $Z_1\in \mathfrak{gl}_{k_1-1}(\mathbb
C),\,\,Z_2\in \mathfrak{gl}_{l_1}(\mathbb
C)$. Then the operator $\psi(S)$, where $S\in P$, see (\ref{P_osp_m-1}), acts in the fiber $W$  in the following way
$$
\begin{array}{c}
\left(
\begin{array}{cc}
(A_{1}^T)^{-1} & 0 \\
0 & (A_{2}^T)^{-1}
\end{array}
\right) \left(                              \left(
\begin{array}{cc}
Z_1 & Q_1  \\
Q_2 & Z_2  \\
\end{array}
\right)+<E_{k_1+l_1-1}>\right)
\left(
\begin{array}{cc}
A_{1}^T & 0 \\
0 & A_{2}^T
\end{array}
\right)=\\
\left(
\begin{array}{cc}
(A_{1}^T)^{-1}Z_1A_{1}^T & (A_{1}^T)^{-1}Q_1A_{2}^T \\
(A_{2}^T)^{-1}Q_1A_{1}^T & (A_{2}^T)^{-1}Z_2A_{2}^T
\end{array}
\right)+<E_{k_1+l_1-1}>,
\end{array}
$$
where $A_{1}\in \GL_{k_1-1}(\mathbb C)$, $A_{2}\in \GL_{l_1}(\mathbb C)$.
We have proven the following lemma.

\medskip

\lem \label{lem predsosp_2k_1-1} {\it Assume that $\mathfrak v(\mathcal S)\simeq \mathfrak{pgl}_{k_1-1|l_1}(\mathbb C)$. Then the representation $\psi$ of
	$P$ in the fiber $W$ over  $o=P$
	is completely reducible and has the following form:
	\begin{equation}
		\begin{split}
			\label{predstavl_v_W_0_osp_2k_1-1}  \psi|R =\left\{
			\begin{array}{l}
				\Ad_{\rho_1} + \Ad_{\rho_2} + \rho_{1}\otimes\rho_{2}^* +
				\rho_{1}\otimes\rho^*_{2} + 1\; \text{for}\; k_1>1,l_1 > 0,\\
				\Ad_{\rho_1}\; \text{for}\; k_1>1,\; l_1 = 0,\\
				\Ad_{\rho_2}\; \text{for}\; k_1 = 1,\; l_1 > 0.\end{array}\right.
		\end{split}
	\end{equation}
Here $\rho_1$ and $\rho_2$ are standard representations of 
	$\GL_{k_1-1}(\mathbb C)$ and $\GL_{l_1}(\mathbb C)$, respectively;  $\Ad_{\rho_1}$ and
	$\Ad_{\rho_2}$  are adjoint representations of
	$\GL_{k_1-1}(\mathbb C)$ and $\GL_{l_1}(\mathbb C)$ in $\mathfrak
	{sl}_{k_1-1}(\mathbb C)$ and $\mathfrak {sl}_{l_1}(\mathbb C)$, respectively, and $1$ is the trivial $1$-dimensional representation.$\Box$
}

\medskip

Everything is ready to prove the following theorem.

\medskip

\t \label{teor osp_2k_1-1} {\it Assume that $r>1$, $m=2k_1-1$ and $n=2l_1$. If
	$\mathcal O_{\mathcal S}(\mathcal S_0)\simeq \mathbb C$ and $\mathfrak v(\mathcal S)\simeq \mathfrak {pgl}_{k_1-1|l_1}(\mathbb C)$, then
	$\Ker \mathcal{P}= \{0\}$.}

\medskip

\noindent{\it Proof.} Let us calculate the vector space $\widetilde{\mathcal{W}}_{0}(\mathcal B_0)$ using Theorem \ref{teor borel}.
The representation $\psi$ of
$P$ is computed in Lemma
\ref{lem predsosp_2k_1-1}.
From (\ref{predstavl_v_W_0_osp_2k_1-1}) it follows that the highest weights of this representation have the following form:
\begin{itemize}
	\item $\mu_{1}-\mu_{k_1-1}$, $\mu_{1}-\lambda_{l_1}$,
	$\lambda_{1}- \mu_{k_1-1}$, $\lambda_{1}- \lambda_{l_1}$,
	$0$ for $k_1>2$, $l_1>1$;
	\item $\mu_{1}-\lambda_{l_1}$,
	$\lambda_{1}- \mu_{1}$, $\lambda_{1}- \lambda_{l_1}$,
	$0$ for $k_1=2$, $l_1>1$;
	\item $\mu_{1}-\mu_{k_1-1}$, $\mu_{1}-\lambda_{1}$,
	$\lambda_{1}- \mu_{k_1-1}$,
	$0$ for $k_1>2$, $l_1=1$;
	
	\item  $\mu_{1}-\lambda_{1}$,
	$\lambda_{1}- \mu_{1}$, $0$ for $k_1=2$, $l_1=1$;
	
	\item $\lambda_{1}- \lambda_{l_1}$ for $k_1=1$, $l_1>1$;

	\item the representation space is equal to $\{0\}$ if $k_1=1$, $l_1=1$. 
\end{itemize}
Note that  a flag supermanifold $\mathbf {IF}_{k|l}$ is defined if $l_1>0$, since otherwise $\mathbf {IF}_{k|l}$ is a usual flag manifold.

The dominant weights of the representation $\psi$ are:
\begin{itemize}
	\item  $0$ if
	$k_1\geq 2$, $l_1\geq 1$;
	
	\item we do not have dominant weights if  $k_1=1$, $l_1\geq 1$.
	
\end{itemize}
Hence,
\begin{align*}
\widetilde{\mathcal{W}}_{0}(\mathcal B_0) \simeq\mathbb C \,\, \,\text{for}
\,\,\, k_1\geq 2,\, l_1 \geq 1;\quad
\widetilde{\mathcal{W}}_{0}(\mathcal B_0) =\{0\} \,\,\,\text{for} \,\,\,  k_1=1, \,l_1\geq 1.
\end{align*}
Here $\mathbb C$ is the trivial one dimensional  $\mathfrak {osp}_{2k_1-1|2l_1}(\mathbb
C)_{\bar 0}$-module.

We denote by $\mathcal {H}$ the stationary Lie sub-supergroup of the point $o$, and by $\mathfrak{p}$ we denote its Lie superalgebra. Above we calculated the represenration $\psi$ of $P=\mathcal {H}_0$ in the fiber $W$. Now we need to compute the represenration $\varphi$ of the Lie 
superalgebra $\mathfrak {p}$ in $W$. Any odd one-parameter subsupergroup 
with odd parameter  $\tau$
of $\mathcal {H}$
acts in the coordinates $Z_{I_1}$ at the point $o$ by the following formula:
$$
\left(
\begin{array}{ccccc}
E & 0 &0 & \tau C_{11} & 0 \\
0 & E &0 & \tau C_{21} & \tau C_{22} \\
0 &0 &1& \tau G_4^T&0\\
-\tau C^T_{22} & 0&0 & E & 0 \\
\tau C^T_{21} & \tau C^T_{11}& G_4 & 0& E \\
\end{array}
\right) \left(
\begin{array}{cc}
0 & 0 \\
E_{k_1-1} & 0 \\
0 & 0 \\
0 & 0 \\
0 & E_{l_1} \\
\end{array}
\right)= \left(
\begin{array}{cc}
0 & 0 \\
E_{k_1-1} & \tau C_{22} \\
0 & 0 \\
0 & 0 \\
\tau C_{11}^T & E_{l_1} \\
\end{array}
\right).
$$
Hence for $Z_{I_2}$ we have:
$$
Z_{I_2}\mapsto \left(
\begin{array}{cc}
E_{k_1-1} & \tau C_{22} \\
\tau C_{11}^T & E_{l_1}  \\
\end{array}
\right) Z_{I_2}.
$$
We see that the representation $\varphi$ in the fiber $W$ coinsides with its restristion on Lie sub-superalgebra $
\mathfrak {gl}_{k_1-1|l_1}(\mathbb C)\subset \mathfrak{p}.
$
Note that the representation of 
$\mathfrak {gl}_{k_1-1|l_1}(\mathbb C)$ in $W\simeq \mathfrak
{pgl}_{k_1-1|l_1}(\mathbb C)$ is equivalent to the adjoint representation of $\mathfrak {gl}_{k_1-1|l_1}(\mathbb C)$ in $\mathfrak
{pgl}_{k_1-1|l_1}(\mathbb C)$.

Now we use an argument similar that was used in \cite{ViGL in JA}. 
Let $\pi : \mathcal W\to \widetilde{\mathcal W}_0 = \mathcal W/\mathcal W_{(1)}$ be the natural map and $\pi_{o} : \mathcal W\to W$ be the composition of $\pi$ and of the evaluation map at the point $o$. The following diagram is commutative 
$$
\begin{CD}
\mathcal W(\mathcal B_0)@>{[X,\,\,\cdot
	\,\,]}>> \mathcal W(\mathcal B_0)
\\
@V{\pi_{o}}VV @V{\pi_{o}}VV\\
W@>{\varphi(X)}>> W
\end{CD},
$$
where $X\in \mathfrak {gl}_{k_1-1|l_1}(\mathbb C)\subset \mathfrak p$. Note that the vector space $\mathcal W(\mathcal B_0)$ is an ideal in $\mathfrak v(\mathcal M)$ and in particular it is invariant with respect to the action of $\mathcal H$.
Consider the image $\pi_{o}(\mathcal W(\mathcal B_0))\subset W$. According the calculation above $\pi_{o}(\mathcal W(\mathcal B_0))$ is trivial or isomorphic to $\mathbb C$. From the commutativity of our diagram it follows that 
$
\pi_{o}(\mathcal W(\mathcal B_0)) \subset \mathfrak
{pgl}_{k_1-1|l_1}(\mathbb C)
$
is invariant with respect to the adjoint representation of $\mathfrak {pgl}_{k_1-1|l_1}(\mathbb C)$. In other words, $\pi_{o}(\mathcal W(\mathcal B_0))$ is an ideal in  $\mathfrak {pgl}_{k_1-1|l_1}(\mathbb C)$. 

Analyzing ideals in  $\mathfrak {pgl}_{k_1-1|l_1}(\mathbb C)$, we see that $\pi_{o}(\mathcal W(\mathcal B_0))$ never coincides with non-trivial ideals. Hence, $\pi_{o}(\mathcal W(\mathcal B_0))=\{0\}$.
In other words, we proved that all sections of  $\pi(\mathcal W(\mathcal B_0))$ are equal to $0$ at the point $o$. Since the locally free sheaf $\widetilde{\mathcal W}_0$ is homogeneous, we get that sections from $\pi(\mathcal W(\mathcal B_0))$ are equal to $0$ at any point. Therefore, we have $\pi(\mathcal W(\mathcal B_0))=\{0\}$ and $\mathcal W(\mathcal B_0)_{(0)} \simeq \mathcal W(\mathcal B_0)_{(1)}.$
From Lemma  \ref{lem kerosp_2m-1}, it follows that $\Ker\mathcal{P}=\{0\}$.$\Box$

\subsection{A construction of an isomorphism $\mathbf {IF}_{2k_1-1,k_1-1|2l_1,l_1} \simeq \mathbf {IF}_{2k_1,k_1|2l_1,l_1}^o$}

The existence of an isomorphism  $\mathbf {IF}_{2k_1-1,k_1-1|2l_1,l_1} \simeq \mathbf {IF}_{2k_1,k_1|2l_1,l_1}^o$, where $\mathbf {IF}_{2k_1,k_1|2l_1,l_1}^o$ is a connected component of $\mathbf {IF}_{2k_1,k_1|2l_1,l_1}$, was noticed in \cite{oniosp} without proof. The goal of this section is to prove this statement.

Recall that the Lie superalgebra $\mathfrak {osp}_{2k_1|2l_1}(\mathbb C)$ is defined by the equation $M^{ST}\Gamma+ \Gamma M =0$, where the matrix $\Gamma$  in the standard basis of  $\mathbb C^{2k_1|2l_1}$ is given by
\begin{equation}
\label{matr_grama_osp chet}
\Gamma=\left(
\begin{array}{cccc}
0&E_{k_1} & 0 & 0  \\
E_{k_1}& 0 & 0 & 0  \\
0 & 0 & 0 & E_{l_1} \\
0 & 0 & -E_{l_1} & 0 \\
\end{array}
\right).
\end{equation}
Therefore explicitly, we have 
\begin{equation}
\label{matr_osp_chet} \mathfrak {osp}_{2k_1|2l_1}(\mathbb C)=
\left\{ \left(
\begin{array}{cccc}
A_{11}&A_{12}&C_{11}& C_{12}\\
A_{21} &-A_{11}^T  & C_{21} & C_{22} \\
-C_{22}^T & -C_{12}^T & B_{11} & B_{12} \\
C_{21}^T &   C_{11}^T & B_{21} & -B_{11}^T \\
\end{array}
\right),
\begin{array}{l}
A_{21}^T=-A_{21},\\A_{12}^T=-A_{12},\\
B_{12}^T=B_{12}, \\ B_{21}^T=B_{21}
\end{array}
\right\}.
\end{equation}
Here $A_{11}$, $B_{11}$ are square complex matrices of size $k_1$ and $l_1$, respectively.

As above $\mathcal B:= \mathbf {IF}_{2k_1-1,k_1-1|2l_1,l_1}$ is an isotropic super-Grassmannians of maximal type and the point   $o\in\mathcal B_{\red}$ is given in the chart corresponding to 
(\ref{okrestnost_osp_2m-1}) by the following equations 
$$
X_1=0,\quad Y_1=0, \quad Z_1=0, \quad
\Xi_1=0, \quad \H_1=0, \quad \mathcal{Z}_1=0.
$$
Consider an isotropic super-Grassmannian of maximal type $ \mathcal B^1:=  \mathbf {IF}_{2k_1,k_1|l_0,l_1}$ and the following coordinate matrix of $\mathcal B^1$  
\begin{equation}
\label{okrestnost_osp_2m} 
Z_{I_1}^1=
\left(
\begin{array}{cc}
Z^1_1 & \mathcal{Z}^1_1 \\
E_{k_1} & 0 \\
\H^1_1 & Y^1_1 \\
0 & E_{l_1} \\
\end{array}
\right),
\end{equation}
where $Z^1_1$ and $Y^1_1$ are matrices with even coordinates of sizes 
$k_1\times k_1$ and $l_1\times l_1$, respectively,
$\mathcal{Z}^1_1$ and $\H^1_1$ are matrices with odd coordinates. Moreover, $(Z^1_1)^T =- Z^1_1$, $(Y^1_1)^T=Y^1_1$ and $\H^1_1 = -(\mathcal{Z}^1_1)^T$. Details about this construction can be found in \cite{oniosp,Vi osp and pi}.

Denote by $o^1\in\mathcal B^1_{\red}$ the point that is given in the chart corresponding to (\ref{okrestnost_osp_2m}) by the following equations 
$$
X^1_1=0, \quad Y^1_1=0 \quad \text{and} \quad \H^1_1=0.
$$ 
Let $\mathcal H^1$ be super-stabilizers  of 
$o^1$. Denote by $\mathfrak{p}^1$ the Lie superalgebras of $\mathcal H^1$. Recall that above we denoted by $\mathfrak{p}$ the Lie superalgebras of $\mathcal H$. Explicitely we have 
$$
\mathfrak{p}=\left\{ \left(
\begin{array}{ccccc}
A_{11}&0&0&C_{11}& 0\\
A_{21} &-A_{11}^T & G_{2} & C_{21} & C_{22} \\
-G_2^T &   0 & 0 & G_4^T & 0 \\
-C_{22}^T & 0 & 0 & B_{11} & 0 \\
C_{21}^T &   C_{11}^T & G_{4}^T & B_{21} & -B_{11}^T \\
\end{array}
\right),
\begin{array}{l}
A_{21}^T=-A_{21},\\
B_{21}^T=B_{21}
\end{array}
\right\},
$$
$$
\mathfrak{p}^1=\left\{ \left(
\begin{array}{cccc}
A_{11} & 0 & C_{11} & 0 \\
A_{21} & -A^T_{11} & C_{21} & C_{22} \\
-C_{22}^T & 0 & B_{11} & 0 \\
C_{21}^T & C_{11}^T & B_{21} & -B_{11}^T \\
\end{array}
\right),\,\,
\begin{array}{l}
A_{21}^T=-A_{21}, \\
B_{21}^T=B_{21}
\end{array}
\right\}.
$$

Now we  need to construct an embedding of the Lie supergroup 
$\OSp_{2k_1-1|2l_1}(\mathbb C)$ into the Lie supergroup  $\OSp_{2k_1|2l_1}(\mathbb
C)$, that will induce the isomorphism of the supermanifolds $\mathcal B$ and $\mathcal B^1$. Let us transform the matrix  $\Gamma$, see. (\ref{matr_grama_osp}), to the following form:
$$
\Gamma'=\left(
\begin{array}{ccc}
E_{t} & 0 & 0 \\
0 & 0 & E_{l_1} \\
0 & -E_{l_1} & 0 \\
\end{array}
\right), \,\,\,\text{where}\,\,t=2k_1-1\,\, \text{or}\,\,
2k_1.
$$

Indeed, we can use the transformation $\Gamma'=S^T\Gamma S$,
where
$$
S=\left(
\begin{array}{cccc}
\frac{1}{\sqrt{2}}E_{k_1-1} & \frac{i}{\sqrt{2}}E_{k_1-1} & 0 & 0 \\
\frac{1}{\sqrt{2}}E_{k_1-1} & \frac{-i}{\sqrt{2}}E_{k_1-1} & 0 & 0 \\
0 & 0 & 1 & 0 \\
0 & 0 & 0 & E_{2l_1} \\
\end{array}
\right)\,\,\,\text{for}\,\,\mathfrak {osp}_{2k_1-1|2l_1}(\mathbb C)
$$
and
$$
S=\left(
\begin{array}{ccc}
\frac{1}{\sqrt{2}}E_{k_1} & \frac{i}{\sqrt{2}}E_{k_1} & 0  \\
\frac{1}{\sqrt{2}}E_{k_1} & \frac{-i}{\sqrt{2}}E_{k_1} & 0 \\
0 & 0 & E_{2l_1} \\
\end{array}
\right)\,\,\,\text{for}\,\,\mathfrak {osp}_{2k_1|2l_1}(\mathbb C).
$$

Let us take $A\in \mathfrak {osp}_{t|2l_1}(\mathbb C)$, where $t=2k_1-1$ or
$2k_1$. Then we have
$$
A^{ST}(S^{-1})^{T}\Gamma'S^{-1}+ (S^{-1})^T\Gamma' S^{-1}A=0.
$$
Clearly the matrices $\Gamma'$ and $S^{-1}$ commute. Therefore,
$$
A^{ST}(S^{-1})^{T}S^{-1}\Gamma'+ \Gamma'(S^{-1})^T S^{-1}A=0 \quad \text{and} \quad A \mapsto (S^{-1})^T S^{-1}A.
$$

Then the Lie superalgebra $\mathfrak {osp}_{t|2l_1}(\mathbb C)$  contains all super-matrices $A$ such that $A^{ST}\Gamma'+\Gamma' A=0$. Therefore the Lie superalgebras $\mathfrak {osp}_{2k_1-1|2l_1}(\mathbb
C)$ and $\mathfrak {osp}_{2k_1|2l_1}(\mathbb C)$ in another basis has the following form, respectively:
\begin{equation}
\label{matr osp 2k_1-1'} \mathfrak {osp}_{2k_1-1|2n}(\mathbb C)=
\left\{ \left(
\begin{array}{ccccc}
A_{21} &-A_{11}^T & G_{2} & C_{21} & C_{22} \\
A_{11}&A_{12}&G_1&C_{11}& C_{12}\\
-G_2^T &   -G_1^T & 0 & G_4^T & -G_3^T \\
-C_{22}^T & -C_{12}^T & G_{3} & B_{11} & B_{12} \\
C_{21}^T &   C_{11}^T & G_{4}^T & B_{21} & -B_{11}^T \\
\end{array}
\right)
\begin{array}{l}
A_{21}^T=-A_{21},\\
A_{12}^T=-A_{12}, \!\!\\
B_{12}^T=B_{12}, \\
B_{21}^T=B_{21}
\end{array}
\right\},
\end{equation}
and 
\begin{equation}
\label{matr osp 2k_1'} \mathfrak {osp}_{2k_1|2n}(\mathbb C)=\left\{
\left(
\begin{array}{cccc}
A_{21} & -A^T_{11} & C_{21} & C_{22} \\
A_{11} & A_{12} & C_{11} & C_{12} \\
-C_{22}^T & -C_{12}^T & B_{11} & B_{12} \\
C_{21}^T & C_{11}^T & B_{21} & -B_{11}^T \\
\end{array}
\right)
\begin{array}{l}
A_{21}^T=-A_{21},\\
A_{12}^T=-A_{12}, \\
B_{12}^T=B_{12}, \\
B_{21}^T=B_{21}
\end{array}
\right\}.
\end{equation}
The Lie superalgebras $\mathfrak{p}$ and $\mathfrak{p}_1$ in the different basis have the following form, respectively: 
$$ 
\mathfrak{p}=\left\{ \left(
\begin{array}{ccccc}
A_{21}&-A_{11}^T &G_{2}&C_{21}& C_{22}\\
A_{11} &0 & 0 & C_{11} & 0 \\
-G_2^T &   0 & 0 & G_4^T & 0 \\
-C_{22}^T & 0 & 0 & B_{11} & 0 \\
C_{21}^T &   C_{11}^T & G_{4}^T & B_{21} & -B_{11}^T \\
\end{array}
\right),
\begin{array}{l}
A_{21}^T=-A_{21},\\
B_{21}^T=B_{21}
\end{array}
\right\}
$$
and
$$
\mathfrak{p}^1=\left\{ \left(
\begin{array}{cccc}
A_{21} & -A^T_{11} & C_{21} & C_{22} \\
A_{11} & 0 & C_{11} & 0 \\
-C_{22}^T & 0 & B_{11} & 0 \\
C_{21}^T & C_{11}^T & B_{21} & -B_{11}^T \\
\end{array}
\right),\,\,
\begin{array}{l}
A_{21}^T=-A_{21}, \\
B_{21}^T=B_{21}
\end{array}
\right\}.
$$
We define the map $j:\OSp_{2k_1-1|2l_1}(\mathbb C)\to
\OSp_{2k_1|2l_1}(\mathbb C)$ by the following formula:
$$
X\mapsto \left(
\begin{array}{cc}
1 & 0 \\
0 & X \\
\end{array}
\right),
$$
where $X$ is a coordinate matrix of $\OSp_{2k_1-1|2l_1}(\mathbb C)$.

Clearly, $\d\!j(\mathfrak{p})\subset \mathfrak{p}^1$ and $j(\mathcal H) \subset \mathcal H_1 $. Therefore we have the map $\bar j: \mathcal G/\mathcal H \to \mathcal G^1/\mathcal H^1$, where $\mathcal G=  \OSp_{2k_1-1|2l_1}(\mathbb C)$ and $\mathcal G^1 =  \OSp_{2k_1|2l_1}(\mathbb C)$. It is well-known that the corresponding map $\bar j_{0}: (\mathcal G/\mathcal H)_{0} \to (\mathcal G^1/\mathcal H^1)_{0}^o$ of underlying spaces is an isomorphism. Here we denote by $(\mathcal G^1/\mathcal H^1)_{0}^o$ the connected component of $(\mathcal G^1/\mathcal H^1)_{0}$. Moreover, $(\d \bar j )_o: \mathfrak g/ \mathfrak p \to \mathfrak g^1/ \mathfrak p^1$ is an isomorphism. Therefore $j$ is an isomorphism too.  The following theorem is proved.

\medskip

\t\label{theor isomorphism of Grassmannians}{\sl Let $k_1\geq 1$ and $l_1\geq 0$. Then there exists an isomorphism  $\mathbf {IF}_{2k_1-1,k_1-1|2l_1,l_1} \simeq \mathbf {IF}_{2k_1,k_1|2l_1,l_1}^o$, where $\mathbf {IF}_{2k_1,k_1|2l_1,l_1}^o$ is a connected component of $\mathbf {IF}_{2k_1,k_1|2l_1,l_1}$. Moreover the following diagramm 
	$$
	\begin{CD}
	\mathfrak{osp}_{2k_1-1|2l_1}(\mathbb C)@>{\d\!j}>>
	\mathfrak{osp}_{2k_1|2l_1}(\mathbb
	C)\\
	@V{}VV @V{}VV\\
	\mathfrak{v}(\mathbf {IF}_{2k_1-1,k_1-1|2l_1,l_1})@>{\d \bar j}>> \mathfrak{v}(\mathbf {IF}_{2k_1,k_1|2l_1,l_1}^o)
	\end{CD}
	$$
is commutative. Here vertical arrows are natural actions of Lie superalgebras on supermanifolds, $\d\!j$ is an inclusion, $\d \bar j$ is an isomorphism. $\Box$

}

\medskip

\subsection{Calculation of $\Im \mathcal{P}$}
Recall that as above $\mathcal{P}: \mathfrak{v}(\mathcal M )\to \mathfrak{v}(\mathcal B) $ and $\mu : \mathfrak{osp}_{2k_1-1|2l_1}(\mathbb
C) \to \mathfrak{v}(\mathcal M)$, where $\mathcal B= \mathbf {IF}_{2k_1-1,k_1-1|2l_1,l_1}$ and $\mathcal M= \mathbf {IF}_{k|l}$ is the isotropic flag supermanifold of maximal type. 
By Theorem \ref{teor super=Grassmannians} and Theorem \ref{theor isomorphism of Grassmannians} the following homomorphisms 
$$
\mathfrak{osp}_{2k_1|2l_1}(\mathbb
C) \to  \mathfrak{v}(\mathbf {IF}_{2k_1,k_1|2l_1,l_1}^o)\quad  \text{and} \quad \d \bar j: \mathfrak{v}(\mathbf {IF}_{2k_1-1,k_1-1|2l_1,l_1})\to  \mathfrak{v}(\mathbf {IF}_{2k_1,k_1|2l_1,l_1}^o)
$$ 
are isomorphisms. Let us identify
 $\mathfrak{v}(\mathbf {IF}_{2k_1-1,k_1-1|2l_1,l_1})$ with   $\mathfrak{osp}_{2k_1|2l_1}(\mathbb C)$ using these isomorphisms. 
Then, $\Im \mathcal{P}\subset \mathfrak{osp}_{2k_1|2l_1}(\mathbb
C)$ and $\Im (\mathcal{P} \circ \mu) = \d j(\mathfrak{osp}_{2k_1-1|2l_1}(\mathbb C)) \subset \mathfrak{osp}_{2k_1|2l_1}(\mathbb C) $. 

\medskip

\t\label{theor Im P=osp}
{\sl Assume that $r>1$, $k_1\geq 1$,
	$l_1\geq 1$ and $k'\ne (k_1-1,\ldots, k_1-1,0,\ldots,0)$.
	Then  $\Im
	\mathcal{P}=\mathfrak{osp}_{2k_1-1|2l_1}(\mathbb C)$.}
\medskip

\noindent{\it Proof.} {\bf Step 1.} Let us prove that $(\Im \mathcal{P})_{\bar 0}=
\mathfrak{osp}_{2k_1-1|2l_1}(\mathbb C)_{\bar 0}$. 
We put $\mu^1: 	\mathfrak{osp}_{2k_1|2l_1}(\mathbb C) \to \mathfrak{v}(\mathcal B^1)$. The action $\mu^1$ induces the action $\mu^1_{0}$
of $\mathfrak{osp}_{2k_1|2l_1}(\mathbb C)_{\bar 0}$ on $\mathcal B^1_{0}$, where $\mathcal B^1_{0}$ is the inderlying space of $\mathcal B^1$,  such that the following diagram is commutative:
\begin{equation}
	\label{theta}
	\begin{CD}
		\mathfrak{osp}_{2k_1|2l_1}(\mathbb C)_{\bar
			0}@>{\mu^1}>> \mathfrak{v}(\mathcal B^1)_{\bar 0}\\
		@| @VV{\theta^1}V\\
		\mathfrak{osp}_{2k_1|2l_1}(\mathbb C)_{\bar 0} @>{
			{\mu}^1_{0}}>> \mathfrak{v}(\mathcal B^1_{0})
	\end{CD}.
\end{equation}
Here $\theta^1$ is the natural map. Assume that $k_1\geq 1$, $l_1\geq 1$. Then 
$\mu^1$ is an isomorphism by Theorem \ref{teor super=Grassmannians}.
It is a classical result, see for example \cite{ADima}, that 
${\mu}^1_{0}$ is also an isomorphism. Therefore  $\theta^1$ is an isomorphism too. 

 Denote by $\mathcal M^1$ the isotropic flag supermanifold of maximal type $\mathbf {IF}_{k^1|l}$, where $k^1=(2k_1, k_1, k_2,\ldots, k_r)$.
We define by $\theta$ and $\widehat{\theta}$ the following natural homomorphisms, respectively:
$$
\begin{array}{l}
\theta: \mathfrak{v}(\mathcal B)_{\bar 0} \to
\mathfrak{v}(\mathcal B_{0}) \quad  \text{and}\quad
\widehat{\theta}: \mathfrak{v}(\mathcal M)_{\bar 0} \to
\mathfrak{v}(\mathcal M_{0}).
\end{array}
$$
Since $\mathfrak{v}(\mathcal B)_{\bar 0}\simeq \mathfrak{osp}_{2k_1|2l_1}(\mathbb C)_{\bar 0}$, see Theorem \ref{teor super=Grassmannians}, and $\mathfrak{v}(\mathcal B_{0}) \simeq \mathfrak{osp}_{2k_1|2l_1}(\mathbb C)_{\bar 0}$, see \cite{ADima}, as above we conclude that $\theta$ is an isomorphism. 

 Consider the commutative diagram as  (\ref{theta}) for $\widehat{\theta}$:
$$
\begin{CD}
\mathfrak{osp}_{2k_1-1|2l_1}(\mathbb C)_{\bar 0}@>>>
\mathfrak{v}(\mathcal M)_{\bar 0}\\
@| @VV{\widehat{\theta}}V\\
\mathfrak{osp}_{2k_1-1|2l_1}(\mathbb C)_{\bar 0} @>>>
\mathfrak{v}(\mathcal M_{0})
\end{CD}.
$$
If $k'\ne (k_1-1,\ldots, k_1-1,0,\ldots,0)$, see \cite{ADima}, we conclude that the homomorphism in the second line is an isomorphism. Hence, the homomorphism  $\widehat{\theta}$ is surjective.
 Consider the following commutative diagram 
 $$
\begin{CD}
\mathfrak{v}(\mathcal M)_{\bar 0}@>{\mathcal{P}}>> \mathfrak{v}(\mathcal B)_{\bar 0}\\
@V{\widehat{\theta}}VV @V{\theta}VV\\
\mathfrak{v}(\mathcal M_{0})@>{\mathcal{P}_{0}}>> \mathfrak{v}(\mathcal B_{0})
\end{CD}.
$$
Here $\mathcal{P}_{0}$ is the projection induced by the projection of the bundle  $\mathcal M_{0}\to \mathcal B_{0}$.  Summing up, we have $\mathcal{P}= \theta^{-1}\circ \mathcal{P}_{0}
\circ \widehat{\theta}$. The homomorphism $\mathcal{P}_{0}$ is clearly injective. 
Hence, $(\Im
\mathcal{P})_{\bar 0} = \mathfrak{osp}_{2k_1-1|2l_1}(\mathbb
C)_{\bar 0}$.

\medskip

{\bf Step 2.} Let us find $(\Im \mathcal{P})_{\bar 1}$. Consider the following commutator of two odd matrices from
$\mathfrak{osp}_{2k_1|2l_1}(\mathbb C)\simeq \mathfrak{v}(\mathbf {IF}_{2k_1-1,k_1-1|2l_1,l_1})$:
\begin{align*}
\begin{array}{l}
\left[\left(
\begin{array}{cccc}
0 & 0 & 0 & 0 \\
0 & 0 & A_1 & A_2 \\
0 & -A_2^T & 0 & 0 \\
0 & A_1^T & 0 & 0 \\
\end{array}
\right)\begin{array}{c}
\\
\\
\\
,
\end{array}
\left(
\begin{array}{cccc}
0 & 0 & B_1 & B_2 \\
0 & 0 & 0 & 0 \\
-B_2^T & 0 & 0 & 0 \\
B_1^T & 0 & 0 & 0 \\
\end{array}
\right)\right] =
\\
\left(
\begin{array}{cccc}
0 & -B_2A_1^T+B_1A_2^T & 0& 0 \\
-A_1B_2^T+A_2B_1^T & 0 & 0 & 0 \\
0 & 0 & 0 & 0 \\
0 & 0 & 0 & 0 \\
\end{array}
\right),
\end{array}
\end{align*}
where $A_1$, $A_2$, $B_1$, $B_2$  are matrices of size  $k_1\times
l_1$. Assume that $B_i=(b^i_{st})$, where $b^i_{st}=0$ if $s>1$. In other words, the matrix $B_i$ contains only one non-trivial line, the first line.

The first matrix in this commutator is contained in  $\mathfrak{osp}_{2k_1-1|2l_1}(\mathbb C)\subset \mathfrak{osp}_{2k_1|2l_1}(\mathbb C)$. Hence, it is contained in $\Im \mathcal{P}$ as a fundamental vector field corresponding to the action of $\mathfrak{osp}_{2k_1-1|2l_1}(\mathbb C)$. Further, the result of the commutator is contained in
$\mathfrak{osp}_{2k_1|2l_1}(\mathbb C)_{\bar 0}$, but it is not contained in $\mathfrak{osp}_{2k_1-1|2l_1}(\mathbb C)_{\bar 0}=(\Im
\mathcal{P})_{\bar 0}$. Therefore, the second matrix in the commutator is not contained in $\Im \mathcal{P}$. Hence, $\Im
\mathcal{P}=\mathfrak{osp}_{2k_1-1|2l_1}(\mathbb C)$.$\Box$

\medskip

The main result of the paper is the following.

\medskip

\t \label{osp_2k_1-1_obshch} {\sl  Let $r>1$, $k_1\geq 1$ and
	$l_1\geq 1$, the conditions (\ref{eq_condition no functions}) hold, $\mathfrak{v}(\mathbf{F}_{k'|l'})\simeq \mathfrak {pgl}_{k_1|l_1}(\mathbb
	C)$ and $k'\ne (k_1-1,\ldots, k_1-1,0,\ldots,0)$. Then 
	$\mathfrak{v}(\mathbf F_{k|l}(\mathfrak{osp}_{2k_1-1|2l_1}(\mathbb
	C)))\simeq \mathfrak{osp}_{2k_1-1|2l_1}(\mathbb C)$.} $\Box$

\noindent
E.~V.: Departamento de Matem{\'a}tica, Instituto de Ci{\^e}ncias Exatas,
Universidade Federal de Minas Gerais,
Av. Ant{\^o}nio Carlos, 6627, Caixa Postal: 702, CEP: 31270-901, Belo Horizonte, 
Minas Gerais, BRAZIL,
email: {\tt VishnyakovaE\symbol{64}googlemail.com}

\end{document}